\documentclass[fleqn,10pt]{wlscirep}
\usepackage[utf8]{inputenc}
\usepackage[T1]{fontenc}
\usepackage{algorithm}
\usepackage{algpseudocode}
\usepackage{comment}

\title{CRAWLING: a Crowdsourcing Algorithm on Wheels for Smart Parking}

\newcommand{\DKL}{\mathcal{D}_{\text{KL}}}
\newcommand{\E}{\mathbb{E}}
\newtheorem{Problem}{Problem}

\author[1,*]{Émiland Garrabé}
\author[1,*]{Giovanni Russo}
\affil[1]{Department of Information and Electrical Engineering and Applied Mathematics, University of Salerno, Fisciano, 84084, Italy}

\affil[*]{\{egarrabe,giovarusso\}@unisa.it}

\begin{abstract}
We present the principled design of CRAWLING: a CRowdsourcing Algorirthm on WheeLs for smart parkING. CRAWLING is an in-car service for the routing of connected cars. Specifically, cars equipped with our service are able to {\em crowdsource} data from third-parties, including other cars, pedestrians, smart sensors and social media, in order to fulfill a given routing task. CRAWLING relies on a solid control-theoretical formulation and the routes it computes are the solution of an optimal control problem where cars maximize a reward capturing environmental conditions while tracking some desired behavior. A key feature of our service is that it allows to consider stochastic behaviors, while taking into account streams of heterogeneous data. We propose a stand-alone, general-purpose, implementation of CRAWLING and we show its effectiveness on a set of scenarios aimed at illustrating all the key features of our service.  Simulations show that, when cars are equipped with CRAWLING, the service effectively orchestrates the vehicles, making them able to react online to road conditions, minimizing their cost functions. The code implementing our service and to replicate the numerical results is made openly available.
\end{abstract}
\begin{document}

\flushbottom
\maketitle
%
%
\thispagestyle{empty}

\section{Introduction}
By $2050$, it is expected\cite{nations2014world} that $66$\% of the World population will live in urban environments, with many of our cities set to become {\em megacities}. Hence, also driven by the {\em Internet of Things} (IoT) revolution and the explosion in the amount of available data, many governments and local authorities have adopted (or are in the process of adopting) concepts based on Smart Cities to tackle the key challenges related to the growing intensive urbanization\cite{bigdata}. As a result, much research and development effort is currently devoted to the design of IoT and data-driven services with the goal of making Smart Cities sustainable\cite{SILVA2018697}, capable of improving the quality of life of their Citizens with subsequent societal, economic and environmental impacts. The design of Intelligent Transportation Systems (ITSs) is central to achieve each of these goals and an exciting opportunity is offered by the possibility of designing, and deploying, services leveraging the increasing pervasiveness of connected cars. \\
In this context, we propose the principled design of a data-driven service based on crowdsourcing for the routing of connected cars, CRAWLING (CRowdsourcing Algorithm on WheeLs for smart parkING). Connected cars equipped with CRAWLING are able to crowdsource the skills/services they need in order to fulfill a given routing task. The service is also designed to allow cars equipped with CRAWLING to gather data from other sources, such as other cars, pedestrians, smart sensors or social media. As a result, vehicles equipped with CRAWLING become an active part of a sharing economy\cite{Sharing_book} ecosystem, where vehicles share, and create, {\em skills} that can be used within the service. A key feature of our service is that it is built around a framework that intrinsically allows to consider stochastic behaviors, 
thus explicitly accounting for users' stochasticity and to capture their privacy requirements. Specifically, CRAWLING exploits a recent data-driven control algorithm\cite{russo2020crowdsourcing,garrabe2021design} that returns a randomized behavior for the car (i.e., a probability function) by solving a sequential decision-making problem\cite{GARRABE_2022b} that formalizes the tracking of a desired behavior expressed in probabilistic terms. That is, our service systematically minimizes a cost function consisting of: (i) a reward capturing road/traffic conditions; (ii) a regularizer that biases the solution towards some target behavior, allowing CRAWLING to keep into account the possible preferences of passengers. \\
While CRAWLING can be tailored towards generic routing problems, for concreteness we focus here on the problem of routing cars to find parking spaces. Smart parking allocation is attracting much research attention\cite{fog_parking, blockchain_privacy, QL_parking,8467499}, as roaming for parking is known to be a major source of pollution, traffic and user stress\cite{roman_overko_2021_6410603,Jemmali:2022aa,FAHIM2021e07050}. In a broader context, lack of coordination in transportation systems is believed to be correlated to an increase in circulating vehicles\cite{coord_uber}. To mitigate this, in the context of parking management, urban sensors can be used to monitor the availability of parking spaces and the key challenge becomes that of designing algorithms that are able to  optimally route the flow of vehicles\cite{dahleh_motion}. Currently, many such algorithms rely on deterministic decision-making\cite{route_planning_review}. For example, methods\cite{multiDepotVRP} for solving multi-depot vehicle routing frame the problem via deterministic optimization, using binary decision variables. On the other hand, the use of stochastic frameworks has been gaining traction in the Intelligent Transportation Systems (ITSs) literature, for example in the context of differential privacy where privacy is achieved by injecting noise in the system\cite{10.1007/11761679_29}, or in driver\cite{driver_intent} or pedestrian\cite{pedestrian_prediction} intent prediction. Efforts to leverage users' data in Smart City contexts are also the focus of much attention\cite{bigdata1,bigdata2}, for example for energy\cite{DDelectricity} or taxi\cite{DDtaxi} demand prediction, or congestion analysis\cite{congestion_data}. While such efforts often focus on modelling\cite{mod} and/or artificial intelligence methods\cite{DDia, ML_subway}, some specialized stochastic optimization frameworks have also arisen, for example in parking assignment\cite{stochastic_parking_assignment}. Finally, data-driven ITSs, and the data-driven analysis of ITSs\cite{minimum_parking} have also gained  traction, although heterogeneous datasets, emerging from different information streams that need to be merged, have been singled out as a challenge\cite{DDITS}.  As an example of one such stream, social media monitoring has been recently used to identify trends related to traffic obstructions\cite{crowd_sensing} and monitor traffic in real time\cite{twitter_traffic}.\\
In summary, the key contributions of the paper can be summarized as follows: after formalizing the smart parking problem as a data-driven control problem, we propose the principled design of CRAWLING and present a stand-alone, general-purpose, implementation of the service. The service provably returns the optimal behavior for cars on which it is equipped. Then, we evaluate its computation times and discuss how these are suitable for routing applications. We also show, via simulations leveraging the microscopic traffic simulator SUMO\cite{SUMO2018}, the effectiveness of CRAWLING in orchestrating the use of parking spaces in a University Campus during rush hours in the morning. Simulations show that CRAWLING effectively directs the connected cars,  confirming that cars equipped with our service improve their average time-to-parking, avoiding unfavourable areas of the map and reacting online to environmental changes (reported, in the scenarios we consider, via social media). \\
The paper is organized as follows. In Section \ref{sec:methods}, we begin by introducing the data-driven control framework exploited by CRAWLING, and subsequently describing its design, the implementation and the scenarios used to assess its performance. Then, in Section \ref{sec:results}, we report the results obtained with the service. We do so by first evaluating the computational efficiency and then discussing the performance obtained on the reference scenarios, which were specifically chosen to validate all the key features of the service. Finally, we give concluding remarks, and discuss the relevance of the results, Section \ref{sec:discussion}.

\subsection{Notation}\label{sec:notation}
We note the set of road links, or lanes in what follows, as $\mathcal{X}$. Time is discrete and, at time step $k$, the connected car occupies link $\mathbf{x}_k\in\mathcal{X}$. The probability of a connected car of going from $\mathbf{x}_{k-1}$ to $\mathbf{x}_k$ is denoted by $\pi(\mathbf{x}_k|\mathbf{x}_{k-1})$. The car has available a number of services (or sources) and the $i$-th service/source is  denoted by $\pi^{(i)}(\mathbf{x}_k|\mathbf{x}_{k-1})$. The target behavior of the car represents the behavior that passengers would like to track in the absence of road (or environmental) disruptions: this is denoted by $p(\mathbf{x}_{k}|\mathbf{x}_{k-1})$. Finally, the reward obtained by the car at time $k$ for being on $\mathbf{x}_{k}$ is denoted by $r_k(\mathbf{x}_k)$. As we shall see, the reward is used in our service to capture road/environmental disruptions that might not be known to the vehicle passengers. Also, we denote by  $\pi(\mathbf{x}_1,\ldots,\mathbf{x}_N|\mathbf{x}_0)$ and $p(\mathbf{x}_1,\ldots,\mathbf{x}_N|\mathbf{x}_0)$ the product of the transition probabilities $\pi(\mathbf{x}_k|\mathbf{x}_{k-1})$ and $p(\mathbf{x}_{k}|\mathbf{x}_{k-1})$, respectively. These joint probability functions clearly depend on the initial condition of the car, $\mathbf{x}_0$, i.e. the starting link/lane of the trip. Finally, we let: (i) $\E_p[h(\mathbf{X})] := \sum h(\mathbf{x})p(\mathbf{x})$ be the expectation of $h(\cdot)$ with the sum taken over the support of the probability function $p(\mathbf{x})$; (ii) the Kullback-Leibler (KL) divergence\cite{KL_51} of $p(\mathbf{x})$ with respect to $q(\mathbf{x})$, with $p(\mathbf{x})$ absolutely continuous with respect to $q(\mathbf{x})$ is  $\mathcal{D}_{\text{KL}}\left(p(\mathbf{x}) \mid\mid q(\mathbf{x}) \right):= \sum p(\mathbf{x}) \log\frac{p(\mathbf{x}) }{q(\mathbf{x}) } $ and measures the proximity of the pair $p(\mathbf{x})$, $q(\mathbf{x})$. We use the Matlab-like notation $k_1:k_2$ to denote the set of integers between the integers $k_1$ and $k_2$ (included). Consequently, $\{w_k\}_{k_1:k_2}$ denotes the set $\{w_{k_1},\ldots,w_{k_2}\}$.


\begin{figure}[h!]
    \centering
    \includegraphics[width=0.75\columnwidth]{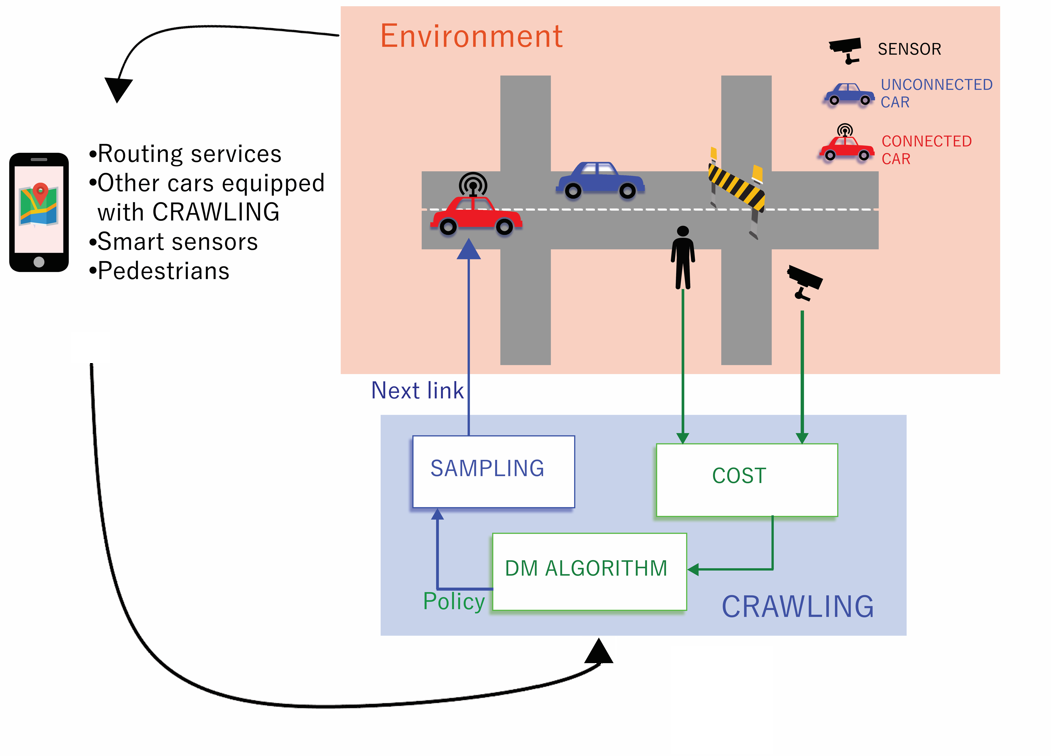}
    \caption{CRAWLING functional architecture and its interactions with the environment. CRAWLING is an in-car service and consists of a decision-making (DM) algorithm, a cost module and a sampling mechanism. These components compute the route the car should take based on environmental data and on the available services.}
    \label{fig:architecture}
\end{figure}

\section{Materials and methods}\label{sec:methods}

We now introduce CRAWLING. We first present the algorithm (Section \ref{sec:service}) and then discuss its implementation (Section \ref{sec:code}) and the reference scenarios used to assess its performance (Section \ref{sec:scenarios}). The key functional components of CRAWLING, and its interactions with the surrounding road environment, are schematically illustrated in Figure \ref{fig:architecture}. As shown in such a figure, CRAWLING is an in-car service: it gathers online data from the environment to feed a suitably defined cost function and has available a number of source services (or simply sources in what follows) to route the car. The sources are e.g., specialized routing services as well as directions collected from other cars in the area also equipped with CRAWLING. Based on the cost and the available sources, CRAWLING dynamically determines the next direction the car should take. This is done via a decision-making (DM) algorithm that exploits recent results in the context of data-driven control\cite{russo2020crowdsourcing,garrabe2021design}. The DM algorithm returns a probability function (i.e., a policy) capturing the turning probabilities of the car. The next road link that the car should take is then obtained by sampling form the policy.

\subsection{CRAWLING: description of the service}\label{sec:service}



CRAWLING seeks to determine the optimal turning probabilities, say $\tilde \pi(\mathbf{x}_k|\mathbf{x}_{k-1})$, of the car by orchestrating the use of a set of available sources. The goal is to track a target behavior for the car, $p(\mathbf{x}_k|\mathbf{x}_{k-1})$, while simultaneously maximizing the car-specific reward $r_k(\mathbf{x}_k)$. The target behavior can be used to capture a preferred set of directions for the car passengers. This can be computed from e.g., historical data\cite{route_pred,8317888} and does not depend on the road conditions. These can be instead captured through the car-specific reward $r_k(\mathbf{x}_k)$. Within the service, the reward signal is used to capture data that are not included in the target behavior: typically, these are related to the state of the network, including hazardous conditions, traffic and roadworks (Scenario $1$ described in Section \ref{sec:scenarios}). Such data can be gathered from a network of urban sensors and the reward can also be fed by interactions with social media (Scenario $2$ in Section \ref{sec:scenarios}). For example, in the simulations from Section \ref{sec:results}, a component of the reward favours road links adjacent to an available parking space, encouraging agents to seek such spaces and park. This structure allows the reward to encompass specifications linked to various tasks, meaning CRAWLING is not per se limited to parking management. To introduce our service, say $S$ the number of available sources, i.e.,  $\pi^{(i)}(\mathbf{x}_k|\mathbf{x}_{k-1})$ with $i=1,\ldots,S$ (see Section \ref{sec:notation} for the definitions). These are assumed to be such that $\DKL\left(\pi^{(i)}(\mathbf{x}_k\mid \mathbf{x}_{k-1})\mid\mid p(\mathbf{x}_k\mid \mathbf{x}_{k-1})\right) < +\infty$. This is a standing assumption for our service that has been shown to be non-restrictive in practice\cite{garrabe2021design}. Given this set-up, the problem tackled by CRAWLING to determine the optimal route for the car can be formalized as follows\cite{russo2020crowdsourcing,garrabe2021design}:
\begin{Problem}\label{prob:problem}
Let $N$ be a positive integer, $\tilde{r}_k(\mathbf{x}_{k-1}) := \E_{\pi(\mathbf{x}_k\mid\mathbf{x}_{k-1})}\left[r_k(\mathbf{X}_k)\right]$ and $\pi(\mathbf{x}_0)$, $p(\mathbf{x}_0)$ be priors on the initial conditions. Find the sequence of weights ${\left\{\boldsymbol{\alpha}_k^\ast\right\}_{1:N}}$, with $\boldsymbol{\alpha}^{\ast}_k:=[\alpha_k^{(1)},\ldots,\alpha_k^{(S)}]$ being an $S$-dimensional vector at time step $k$, solving
\begin{equation}\label{eqn:problem}
    \begin{aligned}
    \underset{\left\{ \boldsymbol{\alpha}_k\right\}_{1:N}}{\text{min}}   &\DKL\left(\pi(\mathbf{x}_1,\ldots,\mathbf{x}_N|\mathbf{x}_0)\pi(\mathbf{x}_0)||p(\mathbf{x}_1,\ldots,\mathbf{x}_N|\mathbf{x}_0)p(\mathbf{x}_0)\right) - \sum_{k=1}^N\E_{\pi(\mathbf{x}_{k-1})}\left[\tilde{r}_k(\mathbf{X}_{k-1})\right]\\
    s.t. & \ \pi(\mathbf{x}_{k}|\mathbf{x}_{k-1}) = \sum_{i\in1:S}\alpha_k^{(i)}\pi^{(i)}(\mathbf{x}_{k}|\mathbf{x}_{k-1}), \ \ \ \forall k \\
    & \sum_{i\in1:S}\alpha_k^{(i)} = 1, \ \ \alpha_k^{(i)} \in \{0, 1\},   \ \ \forall k.
    \end{aligned}
\end{equation}
\end{Problem}
In Problem \ref{prob:problem} the cost function formalizes the fact that the goal of CRAWLING is to track the target behavior while maximizing the car-specific reward. The constraints capture the fact that the solution of the problem is a probability function and, in turn, this probability function is determined by picking one of the sources available to CRAWLING. The problem is tackled\cite{russo2020crowdsourcing,garrabe2021design} by Algorithm \ref{alg:crowd}, which embeds both the {\em cost} and {\em DM algorithm} components of the service in Figure \ref{fig:architecture}. The key macro-steps of the algorithm, which outputs the optimal solution to Problem \ref{prob:problem}, are described next.

The inputs to the algorithm (lines $1-6$) are, besides reward and sources, a time horizon, $N$, and the current state of the agent, $\mathbf{x}_{k-1}$. The target behavior is an optional input parameter. Indeed, if not provided to CRAWLING, this is set to the uniform distribution (lines $7-9$). We note that, when this happens, the first term in the cost of the problem in (\ref{eqn:problem}) becomes an entropic regularizer, also widely used in the literature on reinforcement learning\cite{maxent_levine}. The output of the algorithm (lines $10-11$) are the optimal turning probabilities throughout the time horizon. This is essentially the plan that the car will follow. Following the initialization phase, the for loop computes the \emph{reward-to-go} for the car (this is computed by the {\em cost} module in Figure \ref{fig:architecture}) which represents how \emph{promising} a given direction plan is for the car. The reward-to-go is computed via a backward recursion that starts from the last time step in the time horizon.  Intuitively, the time horizon is a measure of \emph{how far in the future} CRAWLING can look to make a decision (see Section \ref{sec:computation} for the corresponding computational aspect). Specifically, in Algorithm \ref{alg:crowd}, the reward signal $r_k(\cdot)$ received by CRAWLING is updated with information from the future time steps through the signal $\hat{r}_k(\cdot)$, which is obtained via backward recursion. Then, the cost of each source's policy is calculated (lines $17 - 18$) and the source incurring the smallest cost is selected, and used as the agent's policy (i.e., turning probabilities) -- this is returned in line $23$. Finally, the next link that the car needs to take is determined by sampling from the {\em Sampling} component illustrated in Figure \ref{fig:architecture}: this component, given the current state of the car, samples from the next turning probability determined within Algorithm \ref{alg:crowd}.

\begin{algorithm}[!t]
\begin{algorithmic}[1]
\caption{CRAWLING in-car service}\label{alg:crowd}
{
\State {\bf Inputs:}
    \State Current state -- $\mathbf{x}_{k-1}$
	\State Time horizon -- $N$
	\State Reward -- $r_k(\mathbf{x}_k)$
	\State Sources -- $\{\pi^{(i)}(\mathbf{x}_k|\mathbf{x}_{k-1})\}_{1:N}$, with $i\in 1:S$
 	\State Target behavior (optional) -- $\{p(\mathbf{x}_k|\mathbf{x}_{k-1})\}_{1:N}$ 
  	\vspace{0.2cm}
    \If {Target behavior not provided}
    \State $\{p(\mathbf{x}_k|\mathbf{x}_{k-1})\}_{1:N}\gets \text{uniform distribution}$ 
  \EndIf 
	\vspace{0.2cm}
	\State {\bf Output:}
	\State Optimal turning probabilities $\{\tilde\pi(\mathbf{x}_k|\mathbf{x}_{k-1})\}_{1:N}$
	\vspace{0.2cm}
	\State {\bf Initialization:}
	\State Set the vectors $\mathbf{a}^T_{N+1}(\mathbf{x}_N)\gets 0$ and $\tilde{\boldsymbol{\alpha}}_{N+1}\gets 0$
	\vspace{0.2cm}
	\For {$k = N$ to $1$}
	\State 	$\hat{r}_k(\mathbf{x}_k) \gets \mathbf{a}^T_{k+1}(\mathbf{x}_k)\tilde{\boldsymbol{\alpha}}_{k+1}$
	\State $\bar{r}_k(\mathbf{x}_k)\gets {r}_k(\mathbf{x}_k) - \hat{r}_k(\mathbf{x}_k)$
	\For {$i=1:S$} 
	\State $a_k^{(i)}(\mathbf{x}_{k-1})\gets\mathcal{D}_{\text{KL}}\left(\pi^{(i)}(\mathbf{x}_k|\mathbf{x}_{k-1})||p(\mathbf{x}_k|\mathbf{x}_{k-1})\right)+\E_{\pi^{(i)}(\mathbf{x}_k|\mathbf{x}_{k-1})}\left[\bar{r}_k(\mathbf{X}_k)\right]$
	\EndFor
	\State $j_k \gets \underset{i}{\text{argmin }} a_k^{(i)}(\mathbf{x}_{k-1})$
	\State $\mathbf{a}_k(\mathbf{x}_{k-1})\gets [a_k^{(1)}(\mathbf{x}_{k-1}),\ldots,a_k^{(S_k)}(\mathbf{x}_{k-1})]^T$
	\State ${\boldsymbol{\alpha}}_k^{\ast} \gets [\alpha_k^{(1)},\ldots,\alpha_k^{(S_k)}]^T$, where $\alpha_k^{(i)}=0$, $\forall i \ne j_k$ and $\alpha_k^{(j_k)}=1$
	\State $\tilde{\pi}(\mathbf{x}_k|\mathbf{x}_{k-1})\gets \pi^{(j_k)}(\mathbf{x}_k|\mathbf{x}_{k-1})$
	\EndFor
	}
\end{algorithmic}
\end{algorithm}



Before reporting on the implementation of CRAWLING, we also highlight here a benefit in computing (and using as input) probability functions. First, the sources used by CRAWLING are intrinsically stochastic and this can be useful when the agents generating these data want to e.g. guarantee some desired level of differential privacy. This privacy mechanism consists in corrupting the information sent by the source with some (typically) Laplacian noise. Hence, it can be captured with a source having a Laplacian probability function. Moreover, Algorithm \ref{alg:crowd} does not determine the next state of the car but rather a probability function. This aspect can be used to guarantee some level of privacy for drivers using CRAWLING, making it possible to both privately share their policy (which can become a source for other vehicles equipped with CRAWLING) and determine the next direction by sampling from this probability.


\subsection{CRAWLING implementation and validation set-up}\label{sec:code}
\begin{figure}
    \centering
    \includegraphics[width=0.75\linewidth]{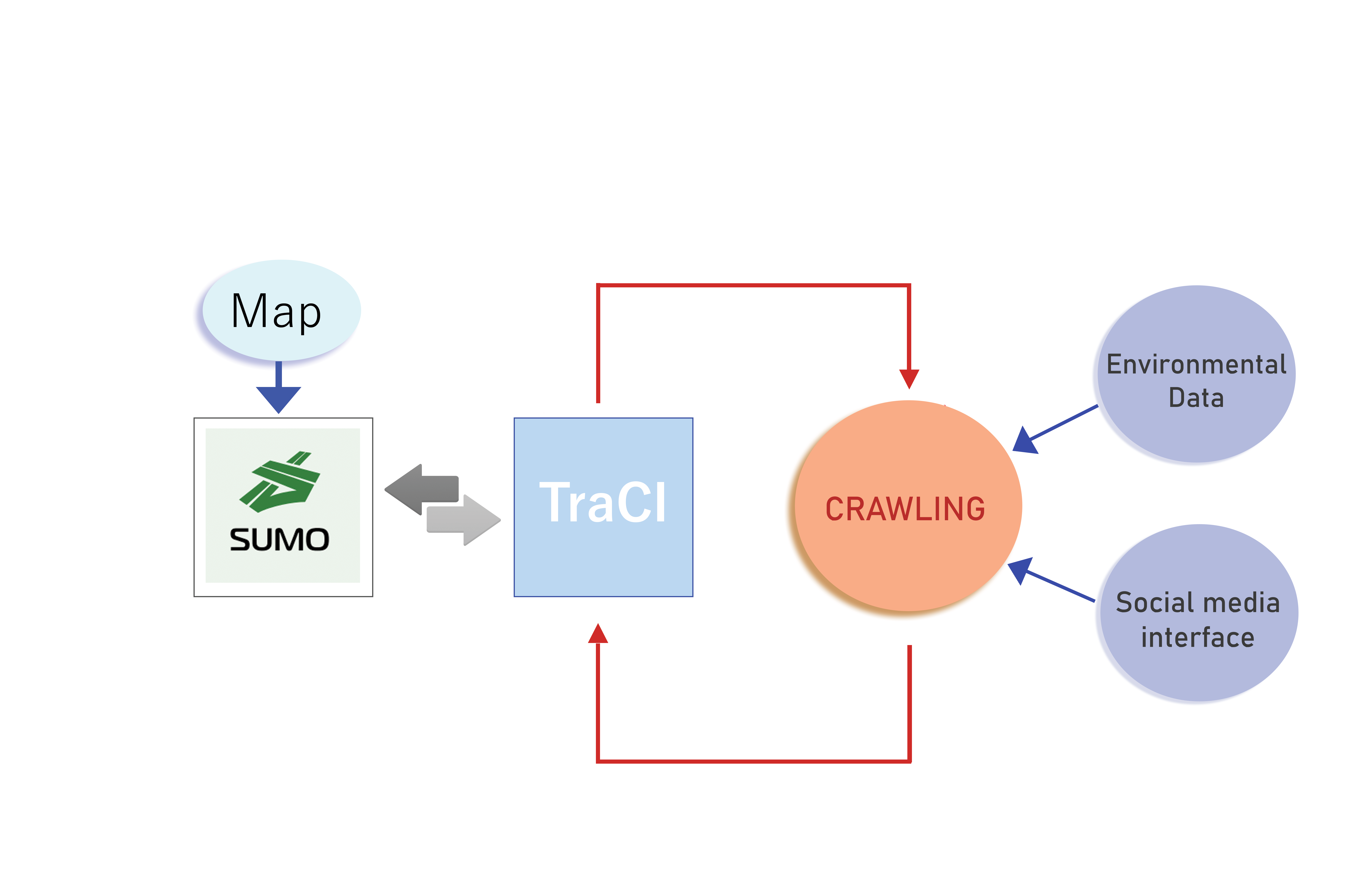}
    \caption{CRAWLING validation set-up. CRAWLING is implemented in Python and interacts with SUMO. This provides the back-end simulation from which CRAWLING gathers the relevant data through TraCI. The directions obtained via Algorithm \ref{alg:crowd} are fed back to the simulation environment via the interface TraCI. For example, in the validation scenarios described in Section \ref{sec:scenarios}, CRAWLING receives data such as parking availability or obstructed road links, which are used to compute directions to connected cars equipped with our service. CRAWLING also interfaces via a Social Interface module with information on social media (in our second scenario, Twitter).}
    \label{fig:setup}
\end{figure}

We now discuss the proposed implementation for CRAWLING. The code can be found at \url{https://tinyurl.com/wcwvcwaj}. The key functional components of the validation set-up are illustrated in Figure \ref{fig:setup}. As shown in such a figure, we validate our service via the microscopic traffic simulator SUMO\cite{SUMO2018}. This serves as environment (in which real maps can be imported and traffic demand can be generated) with which CRAWLING interacts. Interactions with SUMO happen via the Python interface TraCI: this is used to both gather the environmental data from SUMO and to control the vehicles based on the directions provided by Algorithm \ref{alg:crowd}. Moreover, the {\em Social Interface} module gives CRAWLING the possibility of monitoring certain social media (in our scenarios of Section \ref{sec:scenarios}, Twitter). In the set-up we propose, this is done by using the Python library Tweepy (https://www.tweepy.org/), which allows to retrieve a given user's tweet history.

In the CRAWLING implementation, the discrete time steps are not uniform but are rather associated to the vehicle changing link (i.e., when the vehicle transitions state). This design solution removes the need for the cars to have a synchronized clock. Algorithm \ref{alg:crowd} is then implemented via a receding horizon strategy, thus allowing to handle possible changes in the environment. Specifically, every time (say, $t-1$) a given car equipped with CRAWLING transitions to a new link the algorithm: 
\begin{itemize}
    \item gathers the available data and sources (in the setting of Figure \ref{fig:setup} this is done by parsing the data from SUMO through TraCI and from Twitter via Tweepy);
    \item builds the reward based on new information that might have become available during the transition to the new link;
    \item computes the optimal plan of actions for the next $N$ links according to Algorithm \ref{alg:crowd};
    \item determines the next direction by sampling from  $\tilde\pi(\mathbf{x}_{t}|\mathbf{x}_{t-1})$. The car is then controlled through TraCI by setting the direction obtained  from the sampling. 
\end{itemize}

In our code (see our gitub repository for the details) connected cars are encoded by an \emph{agent} Python class, which contains information about the user's goal (i.e., the route/destination the driver vehicle would like to follow/attain). Algorithm \ref{alg:crowd} used by CRAWLING also has a Python class, which takes an agent as argument (i.e., a vehicle). This class regroups all the methods used for the algorithm, namely the computation of the KL divergence, a receding horizon control loop and a sampling mechanism. A simulation file is also provided to manage the key interfaces between the agent, CRAWLING and the simulator used to assess its performance (i.e, SUMO). Specifically, the simulation file includes functions for updating the state of the CRAWLING vehicle in the simulation, updating the reward and performing simulation steps through the SUMO-Python interface TraCI. The simulation file for Scenario $1$ (see Section \ref{sec:scenarios}) also logs results when the simulation is over, while the simulation file for Scenario $2$ (see Section \ref{sec:scenarios}) contains the necessary code to interface CRAWLING with social media, namely Twitter. The file implements a parsing function that, monitoring specific users (one of the authors in this case), is called at every step. If the tweet parsed by our code has the hashtag '\#sumo\_experiment' it is considered for parsing. In this case, the function splits the string into individual words and searches for the word 'blocked' (we leave the design of a more refined parser for future research) and for a road link identifier. Then, the date of the tweet is examined to ensure the tweet has an impact on the current trip: if the tweet was posted before the current date it is ignored. Finally, if all the conditions are met, the blocked road link is assigned a very negative reward as described in Section \ref{sec:sim}.




\begin{figure}
    \centering
    \includegraphics[width=0.75\linewidth]{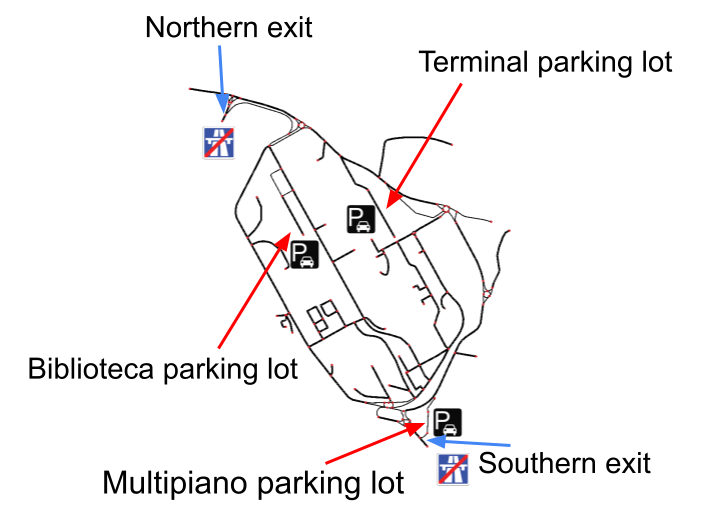}
    \caption{Campus road network. Highways and parking lots ate represented by their respective international symbols.}
    \label{fig:unisa_map}
\end{figure}

\subsection{Scenarios description}\label{sec:scenarios}
We now describe the scenarios (Scenario $1$ and Scenario $2$ in what follows) that are used to validate CRAWLING. The two scenarios were selected to validate different aspects of our service. Namely, with Scenario $1$ we wanted to verify the effectiveness of CRAWLING in smoothly managing a flow of cars that, without our service, could possibly create road blocks. As such, this scenario is centered around the management of a fleet of cars entering the University of Salerno campus during morning rushes. Scenario $2$ was instead aimed at verifying how CRAWLING would adapt to incoming information about road closures made available via a Twitter account.

\noindent \textbf{Scenario $1$.} This scenario simulates a morning rush hour at the campus of the University of Salerno. The campus of the University of Salerno, shown in Figure \ref{fig:unisa_map} is served by two highway exits, and has three parking lots (in our simulations, each parking lot can accommodate up to $50$ cars). The highway exits and corresponding campus entrances will be referred to as southern and northern entrance and the parking lots as Multipiano, Terminal and Biblioteca parking lots (see Figure \ref{fig:unisa_map} for the position of such parking lots on the campus map). In this scenario, a fleet of $150$ cars (i.e., the maximum number of cars that can be accommodated in the simulation) arrives on campus before the teaching begins, with students and faculty staff searching for parking spaces. In particular, $100$ cars arrive from the northern highway exit and seek to park on either Terminal or Biblioteca, while the $50$ remaining cars arrive from the south and reach the Multipiano parking lot. The cars arrive on campus one by one at $15$-second intervals, with the order of arrival being shuffled. That is, when each car arrives, its destination on campus and whether it is equipped with CRAWLING or not, are randomly selected. While cars arrive to campus, roadworks are being carried out on the main ramp leading from the northern highway exit to the campus, rendering traffic on the corresponding links (highlighted in blue in Figure \ref{fig:scenar_osbtru}) difficult (the vehicles using such links have their speed restricted to less than one kilometer per hour). This ramp also corresponds to the preferred route of uncontrolled cars (that is, cars not equipped with CRAWLING) that seek to park in Terminal. In Section \ref{sec:results} we use this scenario to benchmark the effectiveness of CRAWLING in properly managing the flow of cars. Simulations with varying amounts of controlled cars are performed to quantify performance.

\noindent \textbf{Scenario $2$.} In the second scenario, only two cars are considered. They both enter the campus from the northern highway exit, at an interval of $30$ seconds, and seek to reach the Biblioteca parking lot. Both cars are equipped with CRAWLING. Within the scenario, the ramp leading directly from the highway exit to the parking lot is closed. Specifically, after the first car passes through the exit, the ramp is shut off completely and this information is shared through Twitter (by e.g., pedestrians in the area or the local public authority). This information is gathered by CRAWLING and leveraged to re-route the second car entering in the simulation, to find an alternative route that would properly allow it to park in a suitable parking space.


\section{Results}\label{sec:results}
We now describe the results obtained with CRAWLING on the two scenarios of Section \ref{sec:scenarios}. To do so, we first empirically investigate the computational load incurred by  CRAWLING (and discuss our design choices to make the algorithm suitable for routing applications) and then we describe the main settings of our simulation framework. The code can be accessed at \url{https://tinyurl.com/wcwvcwaj}. Note that, to replicate our simulations, a working installation of SUMO, TraCI, Tweepy and the main scientific computing Python libraries is required, along with a Twitter account with developer access. 

\subsection{Computational performance}\label{sec:computation}

\begin{figure}
    \centering
    \includegraphics[width=0.33\columnwidth]{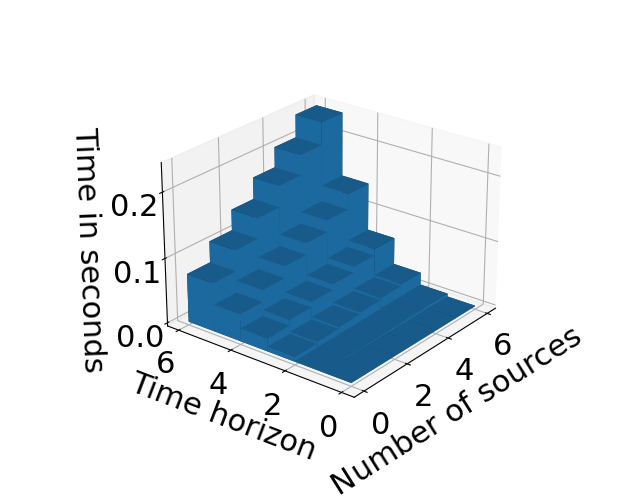}
    \includegraphics[width=0.33\columnwidth]{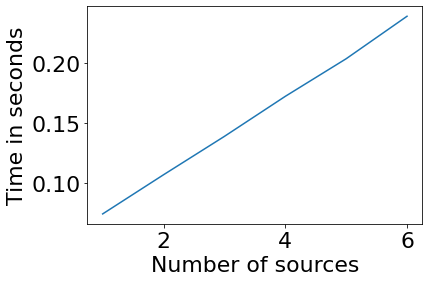}
    \includegraphics[width=0.33\columnwidth]{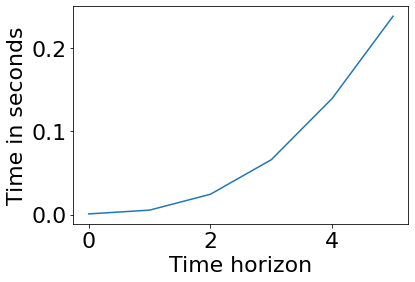}
    \caption{CRAWLING execution time as a function of available sources and time horizon. Leftward panel: computation time (in seconds) as a function of both time horizon and number of services. Middle panel: computation time (in seconds) as a function of the time horizon with $6$ available sources. Rightward panel: computation time (in seconds) as a function of the amount of sources for a time horizon of $5$ steps ahead.}
    \label{fig:comput}
\end{figure}

To make CRAWLING suitable for in-car operation, we leverage the fact that cars being on a given link (i.e., state $\mathbf{x}_{k-1}$) can only transition to an outgoing link (that is, $\mathbf{x}_{k}$ can only be a link in the set of outgoing neighbors of $\mathbf{x}_{k-1}$).  This means that the support of the pfs $\pi(\mathbf{x}_{k}|\mathbf{x}_{k-1})$, $p(\mathbf{x}_{k}|\mathbf{x}_{k-1})$ and $\pi(\mathbf{x}_{k}|\mathbf{x}_{k-1})$ is restricted to the outgoing neighboring links of $\mathbf{x}_{k-1}$. Hence, only the neighbors of $\mathbf{x}_{k-1}$ need to be considered in the computations of Algorithm \ref{alg:crowd}. Thus, in the algorithm implementation: (i) the set of outgoing neighbors associated to the current link of the car is determined; (ii) computations are subsequently performed only over the neighborhood (which depends on the time horizon) rather than over the full state space (which would amount to the full map). Using this approach, we empirically investigated if the service is fast enough for in-car applications.  We did this by running  CRAWLING on each link of the University Campus map in Figure \ref{fig:unisa_map} and logging the average running time. On this map, each road link is connected, on average, to $7$ other links. We repeated these simulations by varying, across them, both the amount of data sources available to the decision-maker (i.e., $S$ was gradually increased between $1$ and $6$) and the time horizon (i.e., $N$ was gradually increased between $0$ and $5$). The results are summarized in Figure \ref{fig:comput}, leftward panel. The figure illustrates that the computation time is linear with respect to the number of sources (Figure \ref{fig:comput}, middle panel). This is due to the fact that Algorithm \ref{alg:crowd} iterates once per source at each time step. More interestingly, the computation time is approximately exponential with respect to the time horizon (Figure \ref{fig:comput}, rightward panel). This is aspect is essentially due to the exponential increase in the state space (even if reduced as described above) as the time horizon, $N$, increases. 

\subsection{Simulation settings}

We now report the settings that are used to implement the scenarios of Section \ref{sec:scenarios} and to obtain the results of Section \ref{sec:sim}.

\noindent \textbf{Scenarios settings.} The parameters are set before running the simulation. In what follows we simply say that vehicles equipped with CRAWLING are {\em controlled} and vehicles without the service are {\em uncontrolled}. First, we created lists containing: (i) the starting points and target parking lots for uncontrolled vehicles; (ii) the target behavior for controlled vehicles. These lists are then shuffled to randomize the departing order of the vehicles in the simulation. For simulations involving both controlled and uncontrolled vehicles, the indexes of controlled and uncontrolled vehicles are randomly assigned. All this information is stored as .npy files. Namely: for Scenario $1$, the files {\em agent.npy} and {\em foe.npy} contain this information. These files are generated before each simulation to avoid biased results. This can be done using the notebook {\em Simulation launcher.ipynb}. The relevant maps were imported from OpenStreetMap (https://www.openstreetmap.org/) using the netconvert software and cleaned up using netedit. A full tutorial for this procedure can be found in the SUMO online documentation. At the beginning of each simulation, SUMO loads the map, while TraCI adds each car from the .npy files to the simulation. The sources, stored as numpy arrays, are loaded when CRAWLING is launched. While SUMO runs the simulation, information including vehicle data and network conditions is obtained in TraCI and transmitted to CRAWLING, which builds the reward.
 
\noindent \textbf{Car settings.} All cars are implemented as TraCI agents. Controlled cars query CRAWLING for directions each time they transition to a new link, while uncontrolled cars are assigned a pre-defined trajectory at the beginning of the simulation (the trajectory is determined by the built-in SUMO routing function that determines the shortest path to the destination). Uncontrolled cars are set to automatically park at the end of their course if possible. On the other hand, controlled cars are assigned a parking space if they arrive to a non-full parking lot. If the parking lot is full, both controlled and uncontrolled cars are rerouted by assigning them a new, non-full parking lot (for uncontrolled cars this is done by assigning them a new destination in SUMO, while controlled agents are assigned a new target behavior).

\begin{figure}
    \centering
    \includegraphics[width=0.75\linewidth]{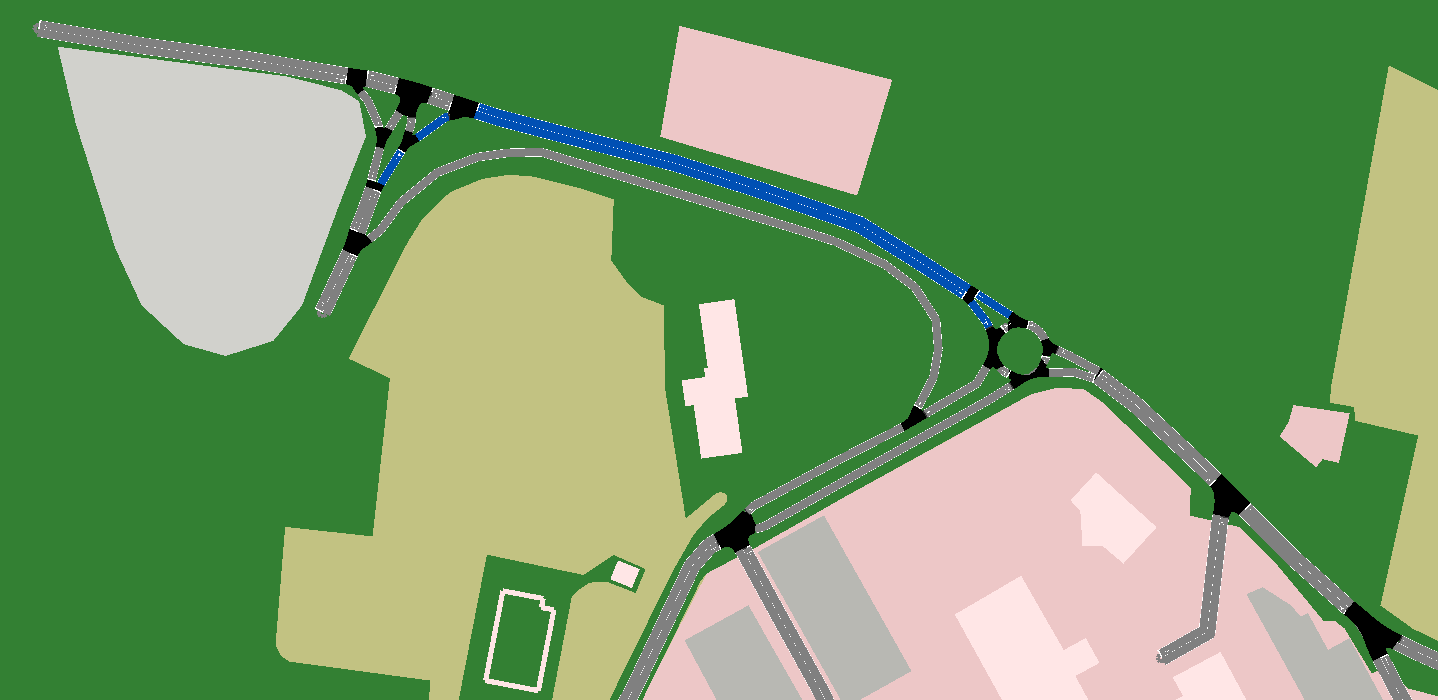}
    \caption{Northern section of the campus with obstructions highlighted in blue for Scenario $1$. Map imported from OpenStreetMap in SUMO via its netconvert tool}
    \label{fig:scenar_osbtru}
\end{figure}

\subsection{Results}\label{sec:sim}

We now describe the simulation specifics and report the results.


\noindent \textbf{Reward and behaviors settings.} The reward encompasses environmental disruptions. In our use-case, this is computed from the state of the road links and is related to the availability of parking spaces. Specifically, at each time step, each link adjacent to a parking lot is assigned a reward of $100$ if the parking lot has vacant spaces, and $-10$ if not. Road links on which traffic is perturbed are assigned an additional negative reward: (i) $-20$ if the traffic is heavily slowed down on the link. This has been simulated, as in Scenario $1$, by enforcing a speed limit of less than one kilometer per hour on the perturbed road links; (ii) $-100$ if the road link is blocked altogether (as in Scenario $2$). Before running the simulations, a set of sources is compiled by using SUMO's routing function: each source is designed to represent directions leading to a different point in the network. We build $6$ sources, using the highway exits, the parking lots and  the north-western side of the campus as destinations. This ensures that the directory of behaviors available to the decision-maker provides a wide coverage of directions from each lane. For each road link, the probability functions $\pi^{(i)}(\mathbf{x}_k\mid\mathbf{x}_{k-1})$ of the sources are built by assigning a high probability to the lane indicated by the routing while the other neighbor links have a small, uniform, probability of being selected. We also built a last source by merging behaviors of cars routed towards the Biblioteca and Terminal parking lot. This source can be thought of as emulating the traces from previous vehicles having navigated the campus. In the following, the target behavior of each connected car is selected among the sources.\\
\noindent \textbf{Scenario $1$ specifics.} Our first set of simulations follows the first scenario outlined in Section \ref{sec:scenarios}. We implemented this scenario by simulating $150$ cars arriving at the campus at $15$-second intervals (with a simulation lasting approximately one hour and twenty minutes of simulated time) and seeking to park. All cars entering the campus via the Southern exit are directed towards the Multipiano parking lot. Uncontrolled cars entering through the Northern exit all seek to park in the Terminal parking lot. On the other hand, the target behavior of controlled cars entering through the Northern exit leads to either Biblioteca or Terminal (this can be interpreted as a target behavior being built from traces of previous CRAWLING cars having similar goals), highlighting the interest of obtaining directions through crowdsourcing stochastic policies. However, the road leading from the northern entrance to the roundabout immediately to its East is obstructed, causing traffic to significantly slow down (as described above, the speed limit on such road links is reduced to less than one kilometer per hour). As specified in Section \ref{sec:scenarios}, this obstruction occurs on the preferred route of uncontrolled cars directed towards Terminal, that is, uncontrolled cars within the $100$ cars entering through the Northern highway exit. This is illustrated in Figure \ref{fig:scenar_osbtru}, which captures a section of the full campus map. Such obstructions can originate from road work, or model a limited capacity of the road link related to e.g. severe climatic conditions. 
In the results shown next, we run different set of simulations with $0$, $50$, $100$ and $150$ controlled cars. This was done to evaluate the {\em penetration} rate of the service, showing its effectiveness as the number of cars equipped with CRAWLING increases.

\begin{figure}
    \centering
    \includegraphics[width=0.75\linewidth]{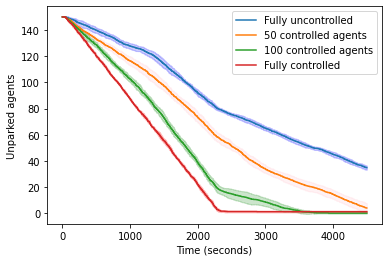}
    \caption{Confidence interval plot of the simulation results. Solid lines represent the mean value, while the shaded area is the standard deviation. Fully controlled fleets manage to fully park within half of the simulation time, while fleets with less controlled cars take longer, failing to find parking spaces for every car before the simulation ends in the extreme case. Note that fleets with just one third of the cars using CRAWLING still significantly outperform the baseline of uncontrolled cars.}
    \label{fig:simres}
\end{figure}

\noindent \textbf{Results for Scenario $1$.} As shown in Figure \ref{fig:simres}, the fully uncontrolled fleet (no vehicles equipped with CRAWLING) fails to park every agent in the allocated time. Indeed, all uncontrolled cars arriving from the northern highway exit are routed to the obstructed lanes. On the other hand, when the fleet is fully controlled (with all vehicles being equipped with our service), CRAWLING only needs approximately half the duration of the simulation to park every car. Namely, CRAWLING allows each controlled cars to adapt its route to avoid the obstruction, effectively orchestrating the use of the available parking spaces. Specifically, we observed that connected cars entering from the Northern exit were routed to the Biblioteca parking lot until it became full. When this happened, CRAWLING rerouted the remaining cars to the Terminal parking lot by a detour route avoiding the obstructed road links. Both controlled and uncontrolled cars entering from the southern highway exit reach the Multipiano parking lot, as it is the closest.
To complement our analysis, we also recorded the average time spent by the average car on an obstructed link and the average time-to-parking over all simulations. These are shown in Table \ref{tab:simres}, together with the results from Welch's t-test\cite{welch_test}, which we used to verify the statistical significance of our simulations. The test is aimed at checking whether two random samples originate from distributions having different means. In other words, a p-value less than $0.05$ would mean that a set of simulations is statistically insignificant and needs to be rejected. Specifically, we use the p-value to compare the average parking time obtained in the five first simulations performed with a set of parameters with the five last simulations following the same parameters.
\begin{table}[]
    \centering
\begin{tabular}[t]{lccc}
\hline
    Scenario & Time (s) spent on obstruction & Average time-to-parking (s) & p-value\\
    \hline
    Fully uncontrolled & 754.2 & 2142.5 & 0.06\\
    50 controlled & 613.1 & 1999.0 & 0.63\\
    100 controlled & 280.2 & 1464.6 & 0.76\\
    150 controlled & 4.1 & 1186.2 & 0.46\\
    \hline
\end{tabular}
    \caption{Summary table for the numerical experiments. All times in seconds.}
    \label{tab:simres}
\end{table}
Table \ref{tab:simres} confirms the qualitative behaviors shown in Figure \ref{fig:simres}. In particular, the increase in the number of cars equipped with CRAWLING decreases both the time spent on obstruction and the time-to-parking. This is due to the fact, with more cars being controlled via CRAWLING, these were able to coordinate to avoid the unfavorable road link, while being directed towards a parking with space with effective capacity to accommodate the incoming cars. 

\begin{figure}
    \centering
    \includegraphics[width=0.75\linewidth]{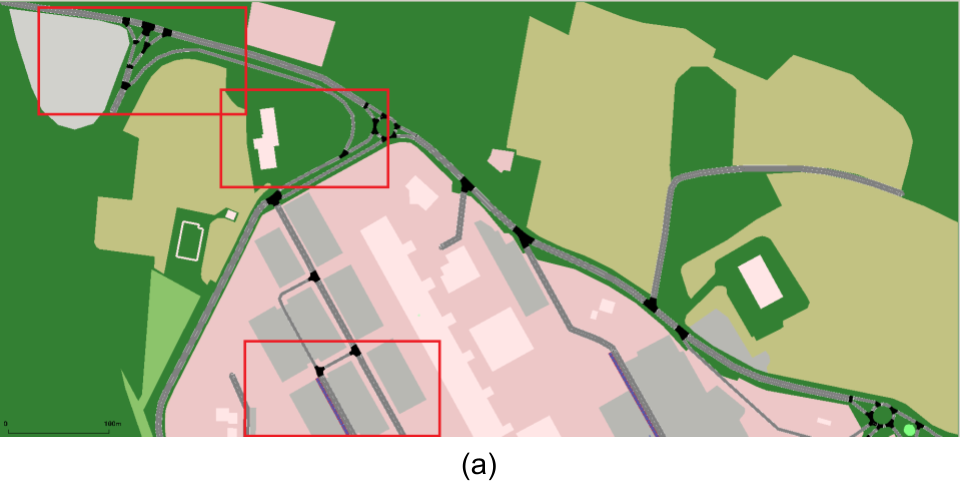}\\
    \includegraphics[width=0.75\linewidth]{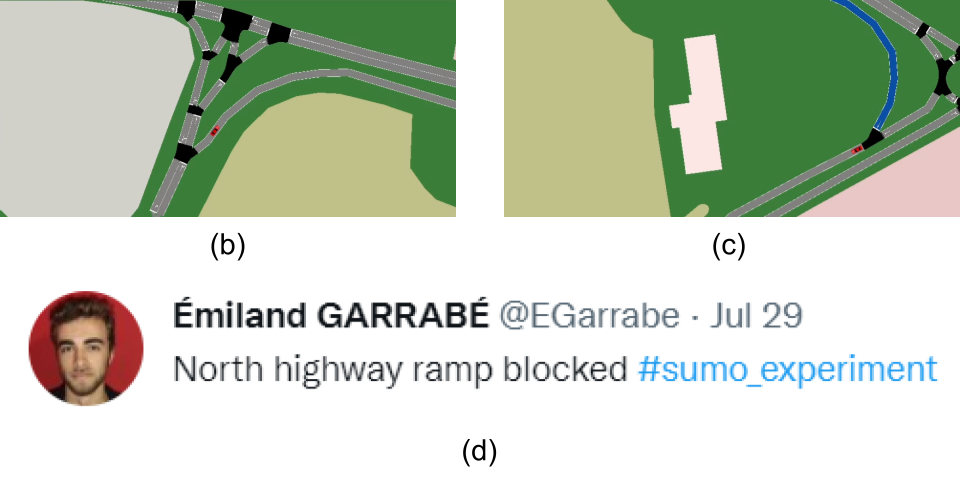}\\
    \includegraphics[width=0.75\linewidth]{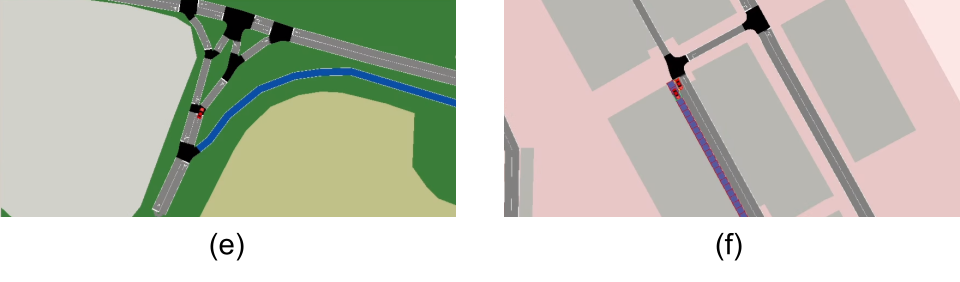}\\
    \caption{Snapshots taken from the simulation for Scenario $2$. Panel (a): areas of the campus interested by the use-case. Specifically, the red rectangles denote the position of the snapshots on the map. Panel (b): first car taking the fastest route from the Northern entrance when there is no obstruction. Panel (c): the route is blocked after the first car has passed the Northern entrance. Panel (d): the information about the disrupting event is shared via a tweet. Panel (e): CRAWLING adapts to handle the disruption and the second car takes an alternate route. Panel (f): both cars reach the desired parking lot.}
    \label{fig:snapshots}
\end{figure}

\noindent \textbf{Scenario $2$ specifics.} For the second scenario in Section \ref{sec:scenarios} we considered two controlled cars entering in the campus from the North entrance, both with the goal of reaching the western parking lot. However, after the first car exits the highway, the road is being blocked. This information, which causes the reward associated to the corresponding road link to sharply decrease, is {\em twitted}. Specifically, from one of the authors' Twitter profile, the following tweet is made: \emph{North highway ramp blocked $\#$sumo\_experiment} after the first car's passage. Periodically, our implementation of CRAWLING checks the Twitter account's post history. The last tweet is parsed and, if an accident is detected -- see Section \ref{sec:code} for a description of the parsing mechanism -- CRAWLING assigns a very negative (i.e., $-100$) reward to the road link where the disrupting event is signalled.\\
\noindent \textbf{Results for Scenario $2$.} The results are summarized in Figure \ref{fig:snapshots}. Such a figure illustrates, in panel (a), the key campus areas for the scenario. At the beginning of the scenario, no road obstruction is detected and indeed, as shown in panel (b), the first car entering in the simulation normally transitions through the Northern exit, following the shortest path through the ramp. After the car passes the ramp, the road obstruction occurs (the obstructed link is highlighted in blue in panel (c)) and the obstruction is reported via a tweet (in this experiment, from one of the authors, see panel (d)). Panel (d) showcases the corresponding tweet that was issued upon detection of the road blocking. Then, panel (e) shows the result of CRAWLING adapting to the new information. Indeed, the route for the second car is re-computed by our service to avoid that the second car enters in the obstructed road link. This is a result of CRAWLING detecting the reward change due to the tweet. An interesting question is to determine a measure of trust for the tweets to be used by CRAWLING. We leave this topic for our future research. At the end of the simulation -- panel (f) -- both cars achieve their goal by reaching their desired parking lot. A video of the corresponding simulation is also available at \url{https://tinyurl.com/wcwvcwaj}.

\section{Discussion and concluding remarks}\label{sec:discussion}
We presented the principled design of CRAWLING, an in-car service for smart parking, for which we make our code available. The service is based on a recent data-driven control algorithm, where the cost function has a Kullback-Leibler divergence component, ensuring that the car tracks a certain target behavior, i.e. a route that the user would like to take. A reward signal is also incorporated into the cost, compiling additional information that is collected through a network of urban sensors and on social media. A key feature of our service is that it determines turning probabilities for cars. Hence, the {\em policy} that is computed by the service is stochastic and it allows to both embed in the framework the possibility that drivers do not follow the indications provided by CRAWLING and to guarantee some desired level of (differential) privacy for users participating to the service. After describing the key functional components of CRAWLING, we proposed a stand-alone, general-purpose implementation of our service and investigated its effectiveness on two Smart Cities scenarios, leveraging SUMO for microscopic traffic simulations. With our first scenario (morning rush) we investigated what happens as the number of vehicles equipped with CRAWLING increases, empirically finding that parking times can be considerably improved even when only $30\%$ of the vehicles were controlled by our service. These simulations were also performed with a limited amount of sources, providing encouraging insights into the system's ability to efficiently make decisions based on a limited set of available services. With our second validation scenario, we investigated the effectiveness of CRAWLING to adapt to  environmental changes reported via social media: simulations once again confirmed the benefits of equipping cars with our service.

In a broader context, a number of considerations can be drawn from our findings. While for concreteness we tailored CRAWLING towards parking management, its approach can be exported to other routing applications. We see our findings as a first step to implement an {\em internet of skills} for connected cars. In this framework, we envisage cars to collaboratively exchange learned skills so as to create a collective intelligence enabling cars to augment each others' abilities through crowdsourcing. We also note that some ventures have already started delivering services based on data collaboratively shared by their users, see for example the Waze application\cite{sharing_economy_mobility}. Interestingly, similar mechanisms based on crowdsourcing can be found in neuroscience -- see, for example the recent {\em thousand brains} theory\cite{1000Brains}, which postulates that the neocortex in our brain does not compute brand new {\em behavioral models} but it patches simpler models from the cortical columns. Finally, in the context of reinforcement learning, a mechanism related to CRAWLING has been recently proposed to orchestrate the use of learning-based and model-based policies to tackle certain control tasks. Interestingly, it has been shown that mechanisms orchestrating these heterogeneous policies can be designed so that they can tutor each other and this, in turn, improves learning efficiency and performance\cite{10.1007/978-3-030-55180-3_6,9654881,pmlr-v168-lellis22a}. Inspired by these studies, possible future directions of our research might involve both embedding tutoring mechanisms in CRAWLING and exploring the links between the approach we presented and certain recent neuroscience theories. We will also embed strategies to filter social media information for their effective use in CRAWLING. Finally, we aim at investigating whether CRAWLING can represent an effective mechanism to fairly regulate possible competitive behaviors between vehicles equipped with our service and legacy, manually-driven, vehicles when these compete for the same, limited, resource.

\section*{Author contributions statement}
G.R. designed the research. E.G. implemented the service, developed the simulations and assessed the numerical methods. Both authors evaluated and discussed the numerical results and contributed in writing the paper. 

\section*{Additional information}

\textbf{Data availability}\\
The code, and the corresponding detailed instructions, needed to replicate the simulations can be found at: \url{https://tinyurl.com/wcwvcwaj}\\
\textbf{Competing interests}\\
The authors do not declare any competing interest.


\bibliography{main}

\end{document}